# Nine Solved and Nine Open Problems in Elementary Geometry

Florentin Smarandache
Math & Science Department
University of New Mexico, Gallup, USA

In this paper we review nine previous proposed and solved problems of elementary *2D* geometry [4] and [6], and we extend them either from triangles to polygons or polyhedrons, or from circles to spheres (from *2D-space* to *3D-space*), and make some comments, conjectures and open questions about them.

**Problem 1.**

We draw the projections $M_i$ of a point $M$ on the sides $A_iA_{i+1}$ of the polygon $A_1...A_n$.
Prove that:

$$\|M_1A_1\|^2 + ... + \|M_nA_n\|^2 = \|M_1A_2\|^2 + ... + \|M_{n-1}A_n\|^2 + \|M_nA_1\|^2$$

**Solution 1.**

For all $i$ we have:

$$\|MM_i\|^2 = \|MA_i\|^2 - \|A_iM_i\|^2 = \|MA_{i+1}\|^2 - \|A_{i+1}M_i\|^2$$

It results that:

$$\|M_iA_i\|^2 - \|M_iA_{i+1}\|^2 = \|MA_i\|^2 - \|MA_{i+1}\|^2$$

From where:

$$\sum_i \left(\|M_iA_i\|^2 - \|M_iA_{i+1}\|^2\right) = \sum_i \left(\|MA_i\|^2 - \|MA_{i+1}\|^2\right) = 0$$

**Open Problem 1.**

1.1. If we consider in a 3D-space the projections $M_i$ of a point $M$ on the *edges* $A_iA_{i+1}$ of a polyhedron $A_1...A_n$ then what kind of relationship (similarly to the above) can we find?

1.2. But if we consider in a 3D-space the projections $M_i$ of a point $M$ on the *faces* $F_i$ of a polyhedron $A_1...A_n$ with k≥4 faces, then what kind of relationship (similarly to the above) can we find?

**Problem 2.**

Let's consider a polygon (which has at least 4 sides) circumscribed to a circle, and $D$ the set of its diagonals and the lines joining the points of contact of two non-adjacent sides. Then $D$ contains at least 3 concurrent lines.

**Solution 2.**

Let $n$ be the number of sides. If $n = 4$, then the two diagonals and the two lines joining the points of contact of two adjacent sides are concurrent (according to Newton's Theorem).

The case $n > 4$ is reduced to the previous case: we consider any polygon $A_1 ... A_n$ (see the figure)

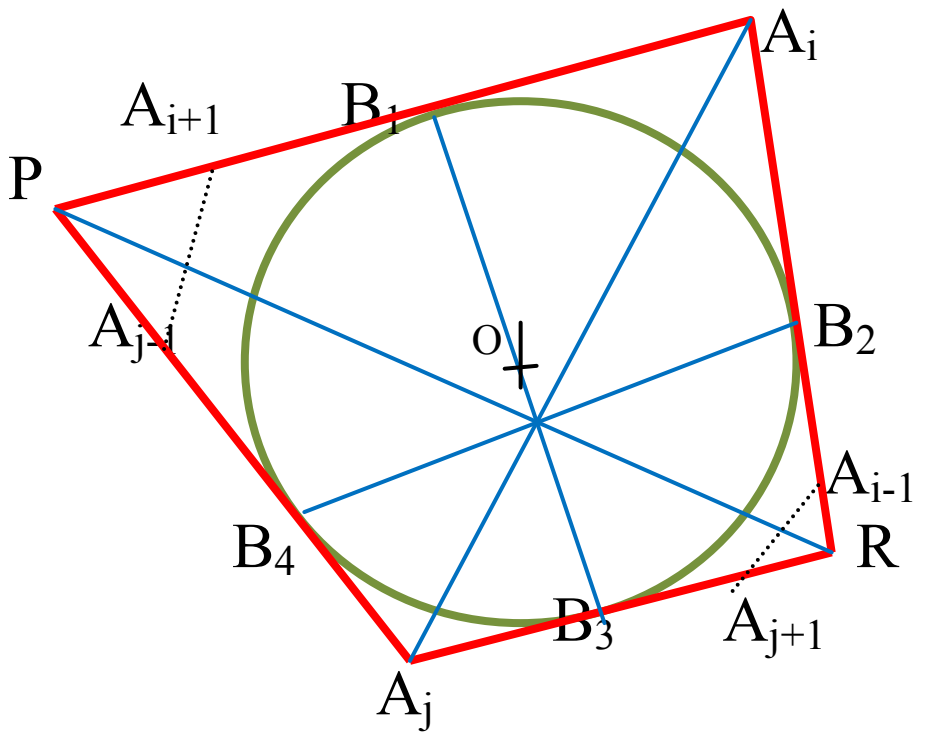

circumscribed to the circle and we choose two vertices $A_i$, $A_j$ $(i \neq j)$ such that

$$A_j\ A_{j-1} \cap A_i\ A_{i+1} = P$$

and

$$A_j\ A_{j+1} \cap A_i\ A_{i-1} = R.$$

Let $B_h$, $h \in \{1, 2, 3, 4\}$ the contact points of the quadrilateral $PA_jRA_i$ with the circle of center $O$. Because of the Newton's theorem, the lines $A_i\ A_j$, $B_1B_3$ and $B_2B_4$ are concurrent.

**Open Problem 2.**

2.1. In what conditions there are more than three concurrent lines?

2.2. What is the maximum number of concurrent lines that can exist (and in what conditions)?

2.3. What about an alternative of this problem: to consider instead of a circle an ellipse, and then a polygon *ellipsoscribed* (let's invent this word, *ellipso-scribed*, meaning a polygon whose all sides are tangent to an ellipse which inside of it): how many concurrent lines we can find among its diagonals and the lines connecting the point of contact of two non-adjacent sides?

2.4. What about generalizing this problem in a 3D-space: a sphere and a polyhedron circumscribed to it?

2.5. Or instead of a sphere to consider an ellipsoid and a polyhedron *ellipsoido-scribed* to it?

Of course, we can go by construction reversely: take a point inside a circle (similarly for an ellipse, a sphere, or ellipsoid), then draw secants passing through this point that intersect the circle (ellipse, sphere, ellipsoid) into two points, and then draw tangents to the circle (or ellipse), or tangent planes to the sphere or ellipsoid) and try to construct a polygon (or polyhedron) from the intersections of the tangent lines (or of tangent planes) if possible.

For example, a regular polygon (or polyhedron) has a higher chance to have more concurrent such lines.

In the 3D space, we may consider, as alternative to this problem, the intersection of planes (instead of lines).

**Problem 3.**
In a triangle $ABC$ let's consider the Cevians $AA'$, $BB'$ and $CC'$ that intersect in $P$. Calculate the minimum value of the expressions:

$$E(P)=\frac{\|PA\|}{\|PA'\|}+\frac{\|PB\|}{\|PB'\|}+\frac{\|PC\|}{\|PC'\|}$$

and

$$F(P)=\frac{\|PA\|}{\|PA'\|}\cdot\frac{\|PB\|}{\|PB'\|}\cdot\frac{\|PC\|}{\|PC'\|}$$

where $A'\in[BC],\ B'\in[CA],\ C'\in[AB]$.

**Solution 3.**
We'll apply the theorem of Van Aubel three times for the triangle $ABC$, and it results:

$$\frac{\|PA\|}{\|PA'\|}=\frac{\|AC'\|}{\|C'B\|}+\frac{\|AB'\|}{\|B'C\|}$$
$$\frac{\|PB\|}{\|PB'\|}=\frac{\|BA'\|}{\|A'C\|}+\frac{\|BC'\|}{\|C'A\|}$$
$$\frac{\|PC\|}{\|PC'\|}=\frac{\|CA'\|}{\|A'B\|}+\frac{\|CB'\|}{\|B'A\|}$$

If we add these three relations and we use the notation

$$\frac{\|AC'\|}{\|C'B\|}=x>0,\ \frac{\|AB'\|}{\|B'C\|}=y>0,\ \frac{\|BA'\|}{\|A'C\|}=z>0$$

then we obtain:

$$E(P)=\left(x+\frac{1}{y}\right)+\left(x+\frac{1}{y}\right)+\left(z+\frac{1}{z}\right)\geq 2+2+2=6$$

The minimum value will be obtained when $x=y=z=1$, therefore when $P$ will be the gravitation center of the triangle.

When we multiply the three relations we obtain

$$F(P) = \left(x + \frac{1}{y}\right) \cdot \left(x + \frac{1}{y}\right) \cdot \left(z + \frac{1}{z}\right) \geq 8$$

**Open Problem 3.**

3.1. Instead of a triangle we may consider a polygon $A_1A_2\ldots A_n$ and the lines $A_1A_1'$, $A_2A_2'$, …, $A_nA_n'$ that intersect in a point P.
Calculate the minimum value of the expressions:

$$E(P) = \frac{\|PA_1\|}{\|PA_1'\|} + \frac{\|PA_2\|}{\|PA_2'\|} + \ldots + \frac{\|PA_n\|}{\|PA_n'\|}$$

$$F(P) = \frac{\|PA_1\|}{\|PA_1'\|} \cdot \frac{\|PA_2\|}{\|PA_2'\|} \cdot \ldots \cdot \frac{\|PA_n\|}{\|PA_n'\|}$$

3.2. Then let's generalize the problem in the 3D space, and consider the polyhedron $A_1A_2\ldots A_n$ and the lines $A_1A_1'$, $A_2A_2'$, …, $A_nA_n'$ that intersect in a point P. Similarly, calculate the minimum of the expressions E(P) and F(P).

**Problem 4.**

If the points $A_1, B_1, C_1$ divide the sides $BC, CA$ respectively $AB$ of a triangle in the same ratio $k > 0$, determine the minimum of the following expression:

$$\|AA_1\|^2 + \|BB_1\|^2 + \|CC_1\|^2 .$$

**Solution 4.**

Suppose $k > 0$ because we work with distances.

$$\|BA_1\| = k\|BC\|, \quad \|CB_1\| = k\|CA\|, \quad \|AC_1\| = k\|AB\|$$

We'll apply tree times Stewart's theorem in the triangle $ABC$, with the segments $AA_1$, $BB_1$, respectively $CC_1$:

$$\|AB\|^2 \cdot \|BC\|(1-k) + \|AC\|^2 \cdot \|BC\|k - \|AA_1\|^2 \cdot \|BC\| = \|BC\|^3(1-k)k$$

where

$$\|AA_1\|^2 = (1-k)\|AB\|^2 + k\|AC\|^2 - (1-k)k\|BC\|^2$$

similarly,

$$\|BB_1\|^2 = (1-k)\|BC\|^2 + k\|BA\|^2 - (1-k)k\|AC\|^2$$

$$\|CC_1\|^2 = (1-k)\|CA\|^2 + k\|CB\|^2 - (1-k)k\|AB\|^2$$

By adding these three equalities we obtain:

$$\|AA_1\|^2 + \|BB_1\|^2 + \|CC_1\|^2 = (k^2 - k + 1)\left(\|AB\|^2 + \|BC\|^2 + \|CA\|^2\right),$$

which takes the minimum value when $k = \frac{1}{2}$, which is the case when the three lines from the enouncement are the medians of the triangle.

The minimum is $\frac{3}{4}\left(\|AB\|^2 + \|BC\|^2 + \|CA\|^2\right)$.

**Open Problem 4.**

4.1. If the points $A_1'$, $A_2'$, …, $A_n'$ divide the sides $A_1A_2$, $A_2A_3$, …, $A_nA_1$ of a polygon in the same ratio k>0, determine the minimum of the expression:

$$\|A_1A_1'\|^2 + \|A_2A_2'\|^2 + ... + \|A_nA_n'\|^2.$$

4.2. Similarly question if the points $A_1'$, $A_2'$, …, $A_n'$ divide the sides $A_1A_2$, $A_2A_3$, …, $A_nA_1$ in the positive ratios $k_1$, $k_2$, …, $k_n$ respectively.

4.3. Generalize this problem for polyhedrons.

**Problem 5.**

In the triangle $ABC$ we draw the lines $AA_1$, $BB_1$, $CC_1$ such that

$$\|A_1B\|^2 + \|B_1C\|^2 + \|C_1A\|^2 = \|AB_1\|^2 + \|BC_1\|^2 + \|CA_1\|^2.$$

In what conditions these three Cevians are concurrent?

**Partial Solution 5.**

They are concurrent for example when $A_1$, $B_1$, $C_1$ are the legs of the medians of the triangle BCA. Or, as Prof. Ion Pătrașcu remarked, when they are the legs of the heights in an acute angle triangle BCA .

More general.

The relation from the problem can be written also as:

$$a\left(\|A_1B\| - \|A_1C\|\right) + b\left(\|B_1C\| - \|C_1A\|\right) + c\left(\|C_1A\| - \|C_1B\|\right) = 0,$$

where $a$, $b$, $c$ are the sides of the triangle.

We'll denote the three above terms as $\alpha$, $\beta$, and respective $\gamma$, such that $\alpha + \beta + \gamma = 0$.

$$\alpha = a\left(\|A_1B\| - \|A_1C\|\right) \Leftrightarrow \frac{\alpha}{a} = \|A_1B\| - \|A_1C\| - 2\|A_1C\|$$

where

$$\frac{\alpha}{a^2} = \frac{a - 2\|A_1C\|}{a} \Leftrightarrow \frac{a^2}{a^2 - \alpha} = \frac{a}{2\|A_1C\|} \Leftrightarrow \frac{a}{2\|A_1C\|} = \frac{2a^2}{a^2 - \alpha} \Leftrightarrow \frac{2a^2 - a^2 + \alpha}{a^2 - \alpha} = \frac{a - \|A_1C\|}{\|A_1C\|}$$

Then

$$\frac{\|A_1B\|}{\|A_1C\|} = \frac{a^2+\alpha}{a^2-\alpha}.$$

Similarly:

$$\frac{\|B_1C\|}{\|B_1A\|} = \frac{b^2+\beta}{b^2-\beta} \text{ and } \frac{\|C_1A\|}{\|C_1B\|} = \frac{c^2+\gamma}{c^2-\gamma}$$

In conformity with Ceva's theorem, the three lines from the problem are concurrent if and only if:

$$\frac{\|A_1B\|}{\|A_1C\|}\cdot\frac{\|B_1C\|}{\|B_1A\|}\cdot\frac{\|C_1A\|}{\|C_1B\|} = 1 \Leftrightarrow \left(a^2+\alpha\right)\left(b^2+\beta\right)\left(c^2+\gamma\right) = \left(a^2-\alpha\right)\left(b^2-\beta\right)\left(c^2-\gamma\right)$$

**Unsolved Problem 5.**

Generalize this problem for a polygon.

**Problem 6.**

In a triangle we draw the Cevians $AA_1$, $BB_1$, $CC_1$ that intersect in $P$. Prove that

$$\frac{PA}{PA_1}\cdot\frac{PB}{PB_1}\cdot\frac{PC}{PC_1} = \frac{AB\cdot BC\cdot CA}{A_1B\cdot B_1C\cdot C_1A}$$

**Solution 6.**

In the triangle $ABC$ we apply the Ceva's theorem:

$$AC_1\cdot BA_1\cdot CB_1 = -AB_1\cdot CA_1\cdot BC_1 \qquad (1)$$

In the triangle $AA_1B$, cut by the transversal $CC_1$, we'll apply the Menelaus' theorem:

$$AC_1\cdot BC\cdot A_1P = AP\cdot A_1C\cdot BC_1 \qquad (2)$$

In the triangle $BB_1C$, cut by the transversal $AA_1$, we apply again the Menelaus' theorem:

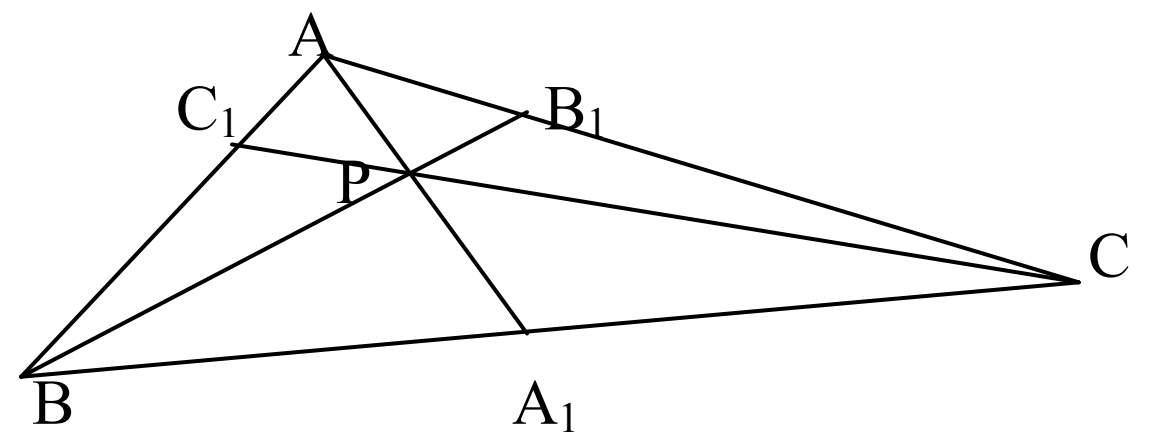

$$BA_1\cdot CA\cdot B_1P = BP\cdot B_1A\cdot CA_1 \qquad (3)$$

We apply one more time the Menelaus' theorem in the triangle $CC_1A$ cut by the transversal $BB_1$:

$$AB\cdot C_1P\cdot CB_1 = AB_1\cdot CP\cdot C_1B \qquad (4)$$

We divide each relation (2), (3), and (4) by relation (1), and we obtain:

$$\frac{PA}{PA_1} = \frac{BC}{BA_1} \cdot \frac{B_1A}{B_1C} \quad (5)$$

$$\frac{PB}{PB_1} = \frac{CA}{CB_1} \cdot \frac{C_1B}{C_1A} \quad (6)$$

$$\frac{PC}{PC_1} = \frac{AB}{AC_1} \cdot \frac{A_1C}{A_1B} \quad (7)$$

Multiplying (5) by (6) and by (7), we have:

$$\frac{PA}{PA_1} \cdot \frac{PB}{PB_1} \cdot \frac{PC}{PC_1} = \frac{AB \cdot BC \cdot CA}{A_1B \cdot B_1C \cdot C_1A} \cdot \frac{AB_1 \cdot BC_1 \cdot CA_1}{A_1B \cdot B_1C \cdot C_1A}$$

but the last fraction is equal to 1 in conformity to Ceva's theorem.

**Unsolved Problem 6.**

Generalize this problem for a polygon.

**Problem 7.**

Given a triangle $ABC$ whose angles are all acute (acute triangle), we consider $A'B'C'$, the triangle formed by the legs of its altitudes.

In which conditions the expression:

$$\|A'B'\| \cdot \|B'C'\| + \|B'C'\| \cdot \|C'A'\| + \|C'A'\| \cdot \|A'B'\|$$

is maximum?

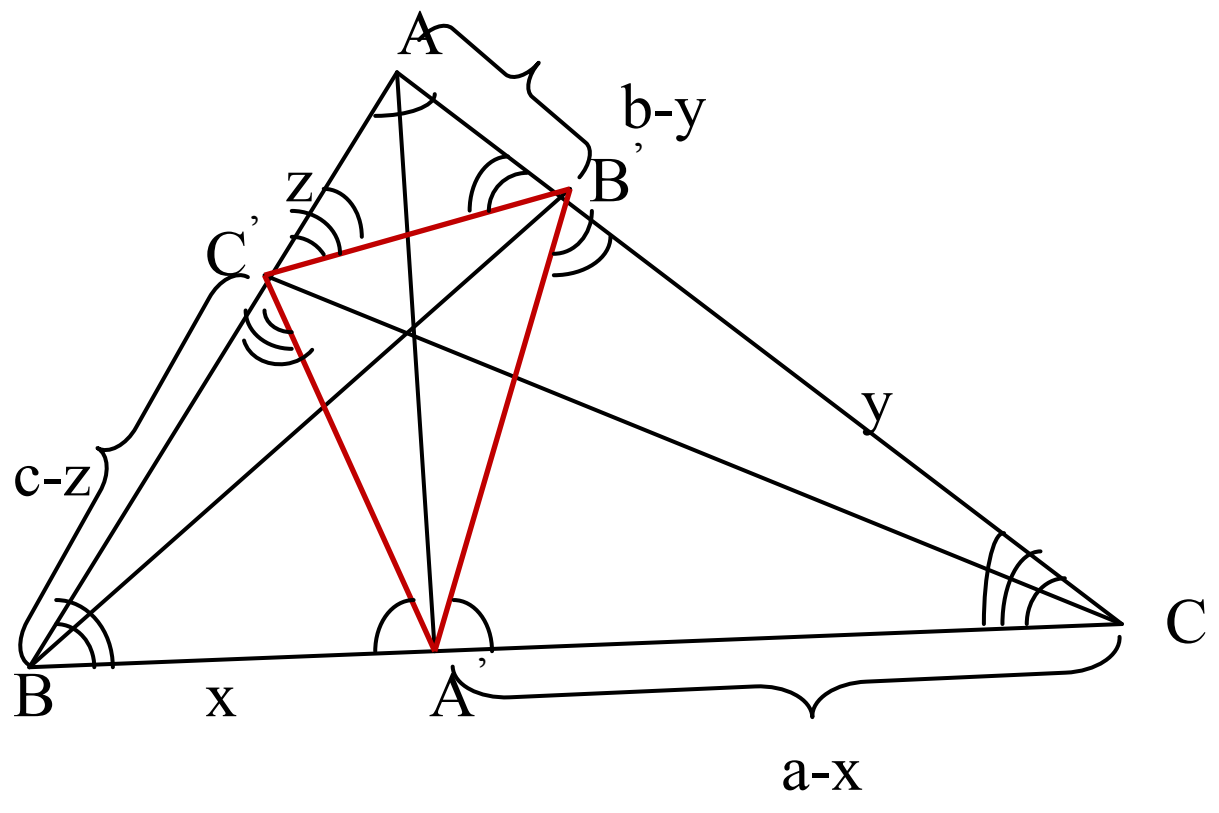

**Solution 7.**

We have

$$\Delta ABC \sim \Delta A'B'C' \sim \Delta AB'C \sim \Delta A'BC' \quad (1)$$

We note

$$\|BA'\| = x, \ \|CB'\| = y, \ \|AC'\| = z .$$

It results that

$$\|A'C\| = a - x,\ \|B'A\| = b - y,\ \|C'B\| = c - z$$

$$\widehat{BAC} = \widehat{B'A'C} = \widehat{BA'C'};\ \widehat{ABC} = \widehat{AB'C'} = \widehat{A'B'C};\ \widehat{BCA} = \widehat{BC'A'} = \widehat{B'C'A}$$

From these equalities it results the relation (1)

$$\Delta A'BC' \sim \Delta A'B'C \Rightarrow \frac{\|A'C'\|}{a-x} = \frac{x}{\|A'B'\|} \quad (2)$$

$$\Delta A'B'C \sim \Delta AB'C' \Rightarrow \frac{\|A'C'\|}{z} = \frac{c-z}{\|B'C'\|} \quad (3)$$

$$\Delta AB'C' \sim \Delta A'B'C \Rightarrow \frac{\|B'C'\|}{y} = \frac{b-y}{\|A'B'\|} \quad (4)$$

From (2), (3) and (4) we observe that the sum of the products from the problem is equal to:

$$x(a-x) + y(b-y) + z(c-z) = \frac{1}{4}\left(a^2 + b^2 + c^2\right) - \left(x - \frac{a}{2}\right)^2 - \left(y - \frac{b}{2}\right)^2 - \left(z - \frac{c}{2}\right)^2$$

which will reach its maximum as long as $x = \frac{a}{2},\ y = \frac{b}{2},\ z = \frac{c}{2}$, that is when the altitudes' legs are in the middle of the sides, therefore when the $\Delta ABC$ is equilateral. The maximum of the expression is $\frac{1}{4}\left(a^2 + b^2 + c^2\right)$.

**Conclusion**[1]: If we note the lengths of the sides of the triangle $\Delta$ABC by ||AB|| = *c*, ||BC|| = *a*, ||CA|| = *b*, and the lengths of the sides of its orthic triangle $\Delta$A\`B\`C\` by ||A\`B\`|| = *c*\`, ||B\`C\`|| = *a*\`, ||C\`A\`|| = *b*\`, then we proved that:

$$4(a`b` + b`c` + c`a`) \le a^2 + b^2 + c^2.$$

**Unsolved Problems 7.**

7.1. Generalize this problem to polygons. Let $A_1A_2\ldots A_m$ be a polygon and P a point inside it. From P, which is called a pedal point, we draw perpendiculars on each side $A_iA_{i+1}$ of the polygon and we note by $A_i$' the intersection between the perpendicular and the side $A_iA_{i+1}$. Let's extend the definition of pedal triangle to a **pedal polygon** in a straight way: i.e. the polygon formed by the orthogonal projections of a pedal point on the sides of the polygon. The pedal polygon $A_1$'$A_2$'…$A_m$' is formed. What properties does this pedal polygon have?

7.2. Generalize this problem to polyhedrons. Let $A_1A_2\ldots A_n$ be a polyhedron and P a point inside it. From P we draw perpendiculars on each edge $A_iA_j$ of the polyhedron and we note by $A_{ij}$' the intersection between the perpendicular and the side $A_iA_{ij}$. Let's name the new formed polyhedron an **edge pedal polyhedron**, $A_1'A_2'\dots A_n'$. What properties does this edge pedal polyhedron have?

---

[1] This is called the **Smarandache's Orthic Theorem** (see [2], [3]).

7.3. Generalize this problem to polyhedrons in a different way. Let $A_1A_2\dots A_n$ be a poliyhedron and P a point inside it. From P we draw perpendiculars on each polyhedron face $F_i$ and we note by $A_i'$ the intersection between the perpendicular and the side $F_i$. Let's call the new formed polyhedron a **face pedal polyhedron**, which is $A_1'A_2'\dots A_p'$, where p is the number of polyhedron's faces. What properties does this face pedal polyhedron have?

**Problem 8.**

Given the distinct points $A_1,\dots,A_n$ on the circumference of a circle with the center in $O$ and of ray $R$.

Prove that there exist two points $A_i$, $A_j$ such that $\left\|\overrightarrow{OA_i}+\overrightarrow{OA_j}\right\| \geq 2R\cos\frac{180^o}{n}$

**Solution 8.**

Because

$$\sphericalangle A_1OA_2 + \sphericalangle A_2OA_3 + \dots + \sphericalangle A_{n-1}OA_n + \sphericalangle A_nOA_1 = 360^o$$

and $\forall i \in \{1,2,\dots,n\}$, $\sphericalangle A_iOA_{i+2} > 0^o$, it result that it exist at least one angle $\sphericalangle A_iOA_j \leq \frac{360^o}{n}$

(otherwise it follows that $S > \frac{360^0}{n}\cdot n = 360^o$).

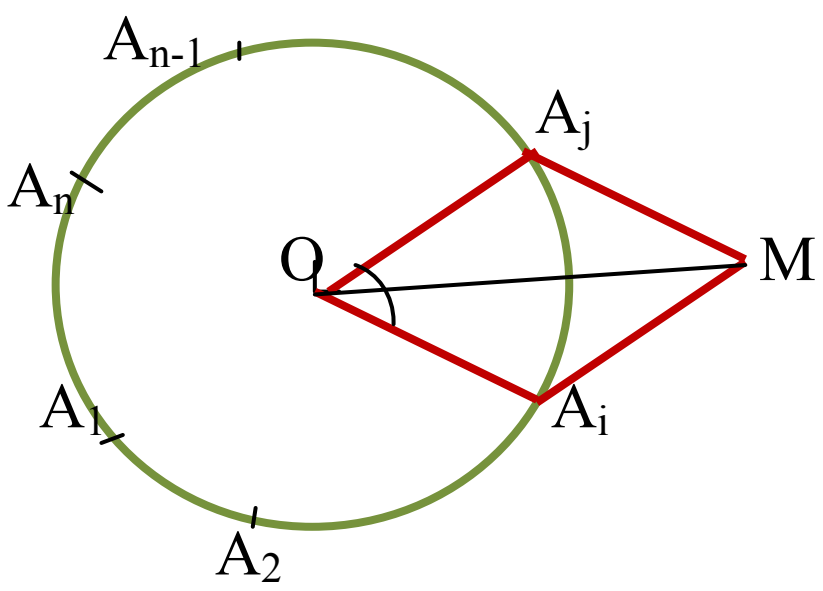

$$\overrightarrow{OA_i}+\overrightarrow{OA_j}=\overrightarrow{OM} \Rightarrow \left\|\overrightarrow{OA_i}+\overrightarrow{OA_j}\right\| = \left\|\overrightarrow{OM}\right\|$$

The quadrilateral $OA_iMA_j$ is a rhombus. When $\alpha$ is smaller, $\left\|\overrightarrow{OM}\right\|$ is greater. As $\alpha \leq \frac{360^o}{n}$, it results that: $\left\|\overrightarrow{OM}\right\| = 2R\cos\frac{\alpha}{2} \geq 2R\cos\frac{180^o}{n}$ .

**Open Problem 8**:

Is it possible to find a similar relationship in an ellipse? (Of course, instead of the circle's radius *R* one should consider the ellipse's axes *a* and *b*.)

**Problem 9:**

Through one of the intersecting points of two circles we draw a line that intersects a second time the circles in the points $M_1$ and $M_2$ respectively. Then the geometric locus of the point $M$ which divides the segment $M_1M_2$ in a ratio $k$ (i.e. $M_1M = k \cdot MM_2$) is the circle of center $O$ (where $O$ is the point that divides the segment of line that connects the two circle centers $O_1$ and respectively $O_2$ into the ratio $k$, i.e. $O_1O = k \cdot OO_2$) and radius $OA$, without the points $A$ and $B$.

**Proof**

Let $O_1E \perp M_1M_2$ and $O_2F \perp M_1M_2$. Let $O \in O_1O_2$ such that $O_1O = k \cdot OO_2$ and $M \in M_1M_2$, where $M_1M = k \cdot MM_2$.

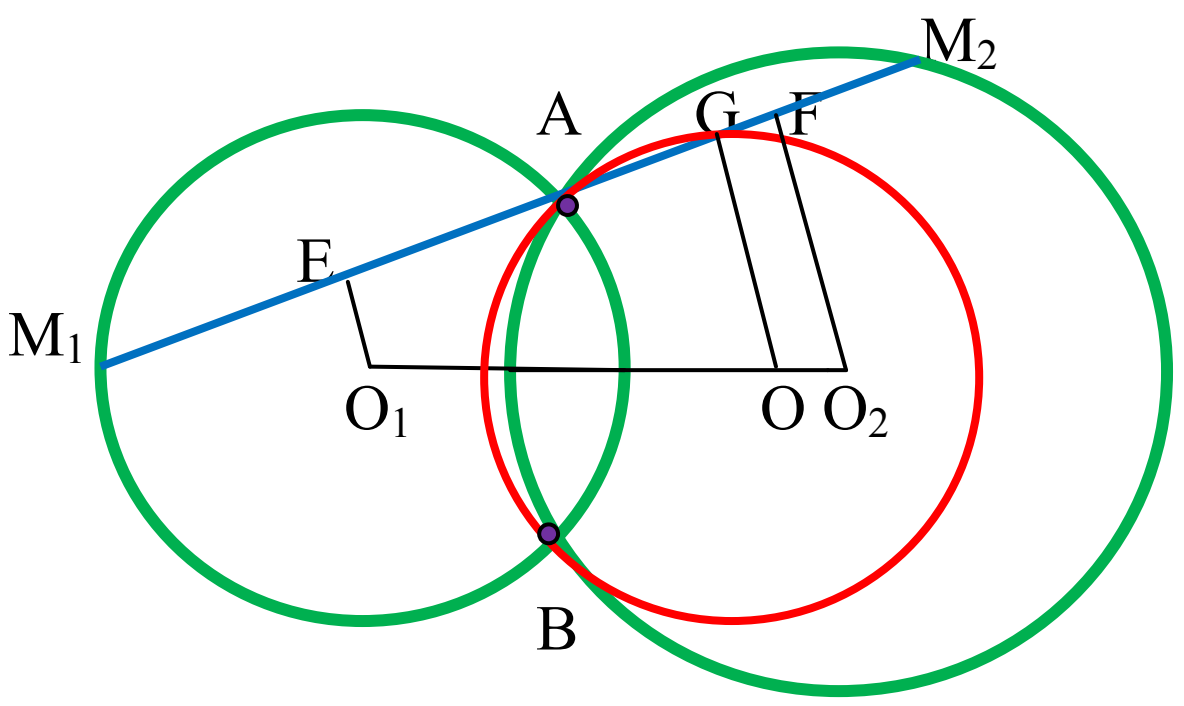

*Fig. 1.*

We construct $OG \perp M_1M_2$ and we make the notations: $M_1E \equiv EA = x$ and $AF \equiv FM_2 = y$.

Then, $AG \equiv GM$, because

$$AG = EG - EA = \frac{k}{k+1}(x+y) - x = \frac{-x+ky}{k+1}$$

and

$$GM = M_1M - M_1A - AG = \frac{k}{k+1}(2x+2y) - 2x - \frac{-x+ky}{k+1} = \frac{-x+ky}{k+1}.$$

Therefore we also have $OM \equiv OA$.

The geometric locus is a circle of center $O$ and radius $OA$, without the points $A$ and $B$ (the red circle in Fig. 1- called *Smarandache's Circle*).

***Conversely.***

If $M \in (GO, OA) \setminus \{A, B\}$, the line $AM$ intersects the two circles in $M_1$ and $M_2$ respectively.

We consider the projections of the points $O_1, O_2, O$ on the line $M_1M_2$ in $E, F, G$ respectively. Because $O_1O = k \cdot OO_2$ it results that $EG = k \cdot GF$.

Making the notations: $M_1E \equiv EA = x$ and $AF \equiv FM_2 = y$ we obtain that

$$M_1M = M_1A + AM = M_1A + 2AG = 2x + 2(EG - EA) =$$
$$= \left[2x + 2\frac{k}{k+1}(x+y) - x\right] = \frac{k}{k+1}(2x+2y) = \frac{k}{k+1}M_1M.$$

For $k = 2$ we find the Problem IV from [5].

**Open Problem 9.**

9.1. The same problem if instead of two circles one considers two ellipses, or one ellipse and one circle.

9.2. The same problem in *3D,* considering instead of two circles two spheres (their surfaces) whose intersection is a circle $\mathcal{C}$. Drawing a line passing through the circumference of $\mathcal{C}$, it will intersect the two spherical surfaces in other two points $M_1$ and respectively $M_2$.

<u>Conjecture</u>: The geometric locus of the point $M$ which divides the segment $M_1M_2$ in a ratio $k$ (i.e. $M_1M = k \cdot MM_2$) includes the spherical surface of center $O$ (where $O$ is the point that divides the segment of line that connects the two sphere centers $O_1$ and respectively $O_2$ into the ratio $k$, i.e. $O_1O = k \cdot OO_2$) and radius $OA$, without the intersection circle $\mathcal{C}$.

A partial proof is this: if the line $M_1M_2$ which intersect the two spheres is the same plane as the line $O_1O_2$ then the *3D* problem is reduce to a *2D* problem and the locus is a circle of radius $OA$ and center $O$ defined as in the original problem, where the point $A$ belongs to the circumference of $\mathcal{C}$ (except two points). If we consider all such cases (infinitely many actually), we get a sphere of radius $OA$ (from which we exclude the intersection circle $\mathcal{C}$) and centered in $O$ *(A* can be any point on the circumference of intersection circle $\mathcal{C}$).

The locus has to be investigated for the case when $M_1M_2$ and $O_1O_2$ are in different planes.

9.3. What about if instead of two spheres we have two ellipsoids, or a sphere and an ellipsoid?